\documentclass[11pt]{article} 
\usepackage{a4,amssymb, amsmath, amsthm}
\usepackage{color}
\usepackage[utf8]{inputenc}
\usepackage[T1]{fontenc}
\usepackage{mathtools}
\usepackage{nicefrac}

 \usepackage{csquotes}
 
 \usepackage{amsmath,amssymb,amsthm,mathrsfs,amsfonts,dsfont}

\newtheorem{theorem}{Theorem}[section]
\newtheorem{corollary}[theorem]{Corollary}

\newtheorem{Definition}[theorem]{Definition}

	\usepackage{todonotes} 
\usepackage{amsmath}
\usepackage{amssymb}
\usepackage{comment}
\usepackage{graphics}
\usepackage{tikz}
\usepackage{pgfplots}
\usepackage[table]{xcolor}
\usepackage{makecell}
\usepackage{graphicx} 
\usepackage{pifont}
\usepackage{algorithm}
\usepackage{algpseudocode}
\usepackage{subcaption}
\usepackage{booktabs}
\usepackage{tikz}
\usepackage{bm,mathrsfs}

\theoremstyle{remark}
\newtheorem{remark}[theorem]{Remark}

\usetikzlibrary{arrows.meta, positioning}

\usepackage{subfig}

\begin{document}

\title{A Stable Boundary Element Method for Reliable Long-Time Industrial Sound Emission}
\author{  Simon Schneider\thanks{Structural Mechanics and Acoustics, Ulm University of Applied Sciences, Germany}$\  {}^\ddag$ \and Ceyhun \"{O}zdemir\thanks{Engineering Mathematics, University of Innsbruck, Austria} \and Heiko Gimperlein${}^\dag$ \and Karsten Urban\thanks{Institute of Numerical Mathematics, University of Ulm, Germany\\ The authors are grateful to A. Aimi and C. Guardasoni for discussions concerning the stability of second-kind
boundary integral formulations.} \and Bernd Graf${}^\ast$ 
}                     
%
%
\date{}

\maketitle \vskip 0.5cm
\begin{abstract}
\noindent In this paper we investigate a stable space-time formulation for long-time industrial sound emission problems. To this end, we use a well-posed Galerkin formulation in space and time of the acoustic wave equation in $\mathbb{R}^3$, involving a hypersingular boundary integral operator. Our numerical experiments confirm that the resulting time stepping scheme is stable and accurate for complex acoustic problems in industrial geometries, in contrast to alternative well-known schemes.
The proposed method is shown to be efficient for real-world problems, and we obtain very good agreement with physical acoustic measurements.  
\end{abstract}

\vskip 1.0cm

\section{Introduction}
\label{sec:intro}

Sound emission and propagation are of central interest in areas such as room acoustics, damage detection, and the prevention of noise pollution. The increasing stringency of noise regulations has increased the impact of Noise, Vibration, and Harshness (NVH) analysis, focused on detecting and mitigating excessive sound radiation caused by structural vibrations, \cite{Vibro,Pang}.
Beyond traditional simulations in the frequency-domain, there is growing interest in time-domain approaches, as they are essential to treat transient situations involving a broad spectrum of frequencies, such as the sequence of acceleration and braking events. 
To predict the sound sources in an NVH analysis in the time-domain, there are stable, accurate and efficient tools, like multibody simulations. However, the calculation of sound radiation in the time-domain has remained largely unresolved, especially in large or unbounded computational domains where the surface area is small compared to the volume. In the frequency-domain, boundary element methods (BEMs) provide established numerical methods. They require computations only on the vibrating surface, thereby reducing the spatial dimension and avoiding the need to construct a mesh of the computational domain. In the time-domain, however, standard BEMs used by engineers are often unstable, limiting their adoption by practitioners.
This work considers a stable BEM for real-world industrial sound emission problems involving complex geometries, long time intervals and a broad spectrum of frequencies. The formulation employs a first-kind boundary integral equation involving the hypersingular operator~$\mathcal{W}$. The weak coercivity of $\mathcal{W}$ assures stability of the associated space-time Galerkin approximations. This article aims to study the performance of this approach in the demanding conditions of real-world engineering problems, including stability, convergence, efficiency and accurate prediction of measurements.\par
\noindent To be specific, we consider the sound radiation of a bounded Lipschitz domain $\Omega\subset \mathbb{R}^3$ resulting from vibrations acting on the boundary $\Gamma:=\partial \Omega$. The radiated sound pressure $u: \mathbb{R}^+ \times \mathbb{R}^3\backslash \overline{\Omega}\rightarrow \mathbb{R}$ satisfies the acoustic wave equation (in the strong form)
\begin{subequations} \label{eq:oriProblem}
	\begin{alignat}{2}
		\frac{1}{c^2}\frac{\partial^2 u}{\partial \tau^2}(\tau,x) -\Delta u(\tau,x)=0, \qquad (\tau,x)\in\mathbb{R}^+ \times \mathbb{R}^3\backslash \overline{\Omega}, \\
		u(\tau,x)=\frac{\partial u}{\partial \tau}(\tau,x)=0, \qquad (\tau,x)\in  \mathbb{R}_0^- \times \mathbb{R}^3\backslash  \overline{\Omega},    \\
		\frac{\partial u}{\partial n}(\tau,x)=f(\tau,x), \hspace{5pt}\quad(\tau,x)\in \mathbb{R}^+ \times \Gamma. \label{eq:ori_problem_NeumannData}
	\end{alignat}
\end{subequations}
Here, $x\in \mathbb{R}^3 \backslash\overline{\Omega}$ or $x\in \Gamma$ is the spatial variable, $\tau\in \mathbb{R}^+$ denotes the time and $c\in\mathbb{R}^+$ the speed of sound. The Neumann data $f$ in \eqref{eq:ori_problem_NeumannData} couple the structural vibrations and the acoustic sound emission:
$$f(\tau,x)=-\rho \frac{\partial^2 s_n}{\partial \tau^2}(\tau,x),$$ 
where $\rho$ is the density of air and $s_n$ the displacement of $\Gamma$ in the outer normal direction, whose second time derivative represents the structure-borne noise.\\  
Problem \eqref{eq:oriProblem} is reduced to a boundary integral equation involving the operator $\mathcal{W}$ (defined in \eqref{eq:hypersingular_equation} below) upon which a space-time Galerkin discretization is performed
{in time and on the two dimensional boundary $\Gamma$.} The chosen discretization leads to a specific structure of the system matrix allowing the space-time problem to be solved by a time-stepping scheme.
Our numerical experiments investigate the sound radiation of two gearbox housings for complex Neumann data obtained from dynamic vibration simulations, considering a harmonic and two transient test cases. The results of the acoustic space-time system are compared with measurements conducted on a test bench in an anechoic chamber. The agreement between the numerical and experimental results highlights the potential of this method for industrial applications.
For comparison and validation, we employ two widely used practical approaches based upon the adjoint double-layer potential operator and a collocation-based method using a combined layer ansatz.\\
In the frequency-domain, BEMs methods have been well-established as an efficient computational method for sound emission problems, whose theoretical and numerical foundations date back to the early 1960s. Initial investigations by Schenck introduced a stabilized approach for practical applications \cite{intro1}, which laid the foundation for many subsequent studies. Nowadays, the frequency-domain BEM is a widely used method in practical engineering applications, as reviewed, for example, in \cite{wu2002boundary,preuss2022recent}.
\\
Our work contributes to the recent growing interests in boundary element methods for wave equations, see \cite{banjai2022integral,costabel04,gimperlein2016adaptive} for an overview of both Galerkin and convolution quadrature methods.
The mathematical analysis of boundary integral formulations of the time-dependent acoustic wave equation was initiated by Bamberger and Ha-Duong \cite{bamberger86}. We refer to \cite{joly2017mathematical} for a refined analysis and to {\cite{aimi2009energy,stut}} for recent developments on the stability of first-kind formulations. 
Key computational advances go back to Terrasse \cite{terrasse93}, in particular the marching-on-in-time (MOT) schemes. Ostermann, Maischak and Stephan analyzed efficient and accurate numerical quadratures for the singular integrals when assembling the system matrix  \cite{gimperlein2016adaptive,Maischak2009,ostermann2010numerical}. 
{Specifically for the hypersingular integral equation, detailed aspects of the quadrature are discussed in \cite{aimi2013neumann}.}
We also highlight recent advances on the efficient assembly and compression of space-time matrices, as well as for higher order methods, applicable both to time-stepping schemes and to more general space-time discretizations \cite{aimi2022partially,aimi2025adaptive,gimperlein2019hp,gimperlein2020residual,polz2021space,veit2016efficient}, and applications to sound emission problems \cite{banz2016time,gimperlein2024space,gimperlein2018time} at low frequencies. {While the numerical stability or instability of time-domain boundary integral equations has been much studied,} the present work considers it in practical applications involving complex Neumann boundary data, complex geometries, and validation against experimental measurements.
\\
\emph{This article is organized as follows:} In Section~\ref{sec:TDBEM}, the Neumann problem of the acoustic wave equation is recalled, followed by the derivation of the first-kind boundary integral equation and the introduction of the corresponding Sobolev function spaces, in which stability and well-posedness can be established. The space-time Galerkin approximation is also presented in this context.
Section~\ref{sec:Algo} describes a suitable discretization through tensor products, outlines the implementation of the hypersingular operator~$\mathcal{W}$, and introduces the MOT algorithm for computing the density.
Section~\ref{sec:ValiStabi} presents numerical experiments on test geometries, including convergence studies for constant space-time ratios and long-time stability investigations using a fixed spatial grid with successively refined time steps.
Section~\ref{sec:Exp} introduces the experimental setup as well as the workflow of a numerical NVH analysis, which serves as the basis for the considered engineering applications presented in Sections~\ref{sec:6} and \ref{sec:7}, where we compare the numerical results of harmonic and transient test cases with experimental measurements. The paper ends with some conclusions in Section~\ref{sec:8}. In the appendix, we collect details of the collocation method and the adjoint double-layer potential approach, with which we compare our method. 
\section{Time-Domain Boundary Integral Formulation for Sound Radiation} \label{sec:TDBEM}
In this section we recall the reformulation of Problem \eqref{eq:oriProblem} for the acoustic wave equation  on the bounded Lipschitz domain $\Omega$ as an equivalent boundary integral equation on the boundary $\Gamma$, following the fundamental work of Bamberger and Ha-Duong \cite{bamberger86}. Well-posedness of the resulting weak formulation and the stability of the boundary element discretization are discussed. 
To simplify notation, we first rescale the time variable \( t = c \, \tau \) in Problem \eqref{eq:oriProblem}:
\begin{align} \label{eq:oriProblem2}
    \frac{\partial^2 u}{\partial t^2}(t,x) - \Delta u(t,x) &= 0, \qquad (t,x)\in \mathbb{R}^+ \times \mathbb{R}^3 \setminus \overline{\Omega},\nonumber\\
    u(t,x) = \frac{\partial u}{\partial t}(t,x) &= 0, \qquad (t,x)\in \mathbb{R}_0^- \times \mathbb{R}^3 \setminus \overline{\Omega}, \\
    \frac{\partial u}{\partial n}(t,x) &= f(t,x),\hspace{3pt} \qquad (t,x)\in \mathbb{R}^+ \times \Gamma\nonumber.
\end{align}
\subsection{Variational Formulation}
We start from a double-layer potential ansatz for the sound pressure \( u \)  in the exterior domain \( \mathbb{R}^+ \times \mathbb{R}^3 \setminus \overline{\Omega} \),  in terms of an unknown acoustic density \( \psi : \mathbb{R}^+ \times \Gamma \rightarrow \mathbb{R}\):
\begin{align}\label{eq:Repr}
  u(t,x) =  D\psi(t,x):=\int_{\mathbb{R}^+}\int_{\Gamma}\frac{\partial G}{\partial n_y}(t-t^*,x,y)\psi(t^*,y) ds_y\,dt^*.
\end{align}
Here, $t^*$ is the time integration variable, $ds_y$ denotes the surface measure on $\Gamma$, and $G:\mathbb{R}\times \mathbb{R}^3\times \mathbb{R}^3\rightarrow \mathbb{R}$ is the fundamental solution
$$G(t, x, y)=\frac{\delta\big(t - |x - y|\big)}{4 \pi |x - y|}$$
of the wave equation in \(\mathbb{R}^3\).
By applying the trace theorem to the normal derivative of the representation \eqref{eq:Repr}, we obtain the first-kind boundary integral equation:
\begin{align} \label{eq:BI1}
\mathcal{W} \psi(t,x) = \frac{\partial}{\partial n} u(t,x)=f(t,x) \quad \text{on } \mathbb{R}^+ \times \Gamma,
\end{align}
which involves the hypersingular operator \(\mathcal{W}\):
\begin{equation}\label{eq:hypersingular_equation}
	\mathcal{W} \psi (t,x):=\int_{\mathbb{R}^+}\int_{\Gamma} \frac{\partial^2 G}{\partial n_x \partial n_y}(t - t^*, x, y)\, \psi(t^*, y)\, dt^*\, ds_y.
\end{equation}
For the analysis of Problem \eqref{eq:BI1}, we introduce space-time anisotropic Sobolev spaces  on the boundary \(\Gamma \subset \mathbb{R}^3\) \cite{gimperlein2017priori,Ha-Duong03a}. 
For screen geometries, where $\partial \Gamma~\neq~\emptyset$, we refer to \cite{gimperlein2017priori} and extend \(\Gamma\) to a closed, orientable Lipschitz manifold \(\widetilde{\Gamma}\), on which the standard Sobolev spaces can be defined:
$$
\widetilde{H}^s(\Gamma)=\{u\in H^s(\widetilde{\Gamma}):\,\text{supp}\,u\subset\bar{\Gamma} \},\quad s\in \mathbb{R}.
$$
In this context, $H^s(\Gamma)$ is the quotient space $H^s(\widetilde{\Gamma})\backslash \widetilde{H}^s(\widetilde{\Gamma}\backslash \bar{\Gamma}).$
To define a specific family of Sobolev norms, we introduce a partition of unity \( \{\alpha_i\}_{i=1,...,p} \) associated with an open covering \( \{B_i\}_{i=1,...,p} \) of the extended surface \( \widetilde{\Gamma} \).  
Each set \( B_i \) is mapped by a diffeomorphism \( \zeta_i \) onto the unit square as a subset of \( \mathbb{R}^{3} \).  
This allows us to define Sobolev norms on \( \widetilde{\Gamma} \) by expressing the localized functions \( \alpha_i u \) in Euclidean coordinates:
\[
\|u\|_{s, \omega, \widetilde{\Gamma}} = \left( \sum_{i=1}^p \int_{\mathbb{R}^d} \left(|\omega|^2 + |\xi|^2\right)^s \left| \mathcal{F} \left\{ (\alpha_i u) \circ \zeta_i^{-1} \right\}(\xi) \right|^2 \, d\xi \right)^{1/2}.
\]
Here, $\mathcal{F}$ denotes the Fourier-Laplace transform. Note that these norms are equivalent for different parameter $\omega \in \mathbb{C}\backslash\{0\}$. 
They induce norms on \( H^s(\Gamma) \) and $ \widetilde{H}^s(\Gamma) $:
\[
\|u\|_{s,\omega,\Gamma} = \inf_{v \in \widetilde{H}^s(\widetilde{\Gamma} \setminus \Gamma)} \|u + v\|_{s,\omega,\widetilde{\Gamma}}, \quad
\|u\|_{s,\omega,\Gamma,*} = \|e^+ u\|_{s,\omega,\widetilde{\Gamma}},
\]
where \( e^+ u \)  denotes the extension of \( u \)  by zero from $\Gamma$ to $\widetilde{\Gamma}$.
The space-time anisotropic Sobolev spaces are now defined as follows:
\begin{Definition}
For $\sigma > 0$ and $r, s \in \mathbb{R}$, define
$$H^r_\sigma(\mathbb{R}^+, H^s(\Gamma)) = \left\{ u \in \mathcal{D}'_+(H^s(\Gamma)) : e^{-\sigma t} u \in \mathcal{S}'_+(H^s(\Gamma)),\, \|u\|_{r,s,\sigma,\Gamma} < \infty \right\},$$
$$H^r_\sigma(\mathbb{R}^+, \widetilde{H}^s(\Gamma)) = \left\{ u \in \mathcal{D}'_+(\widetilde{H}^s(\Gamma)) : e^{-\sigma t} u \in \mathcal{S}'_+(\widetilde{H}^s(\Gamma)),\,\|u\|_{r,s,\sigma,\Gamma,*} < \infty \right\}, $$
where $\mathcal{D}'_+(E)$ and $\mathcal{S}'_+(E)$ denote the spaces of distributions and tempered distributions on $\mathbb{R}$ with support in $[0,\infty)$, taking values in a Hilbert space $E$, respectively.\\ For $E=H^s(\Gamma)$ and $E=\widetilde{H}^s(\Gamma)$, we get the norms:
	$$
	\|u\|_{r,s,\sigma} := \|u\|_{r,s,\sigma,\Gamma} = \left( \int_{-\infty + i\sigma}^{+\infty + i\sigma} |\omega|^{2r} \|\hat{u}(\omega)\|_{s,\omega,\Gamma}^2 \, d\omega \right)^{1/2},
	$$
	$$
	\|u\|_{r,s,\sigma,*} := \|u\|_{r,s,\sigma,\Gamma,*} = \left( \int_{-\infty + i\sigma}^{+\infty + i\sigma} |\omega|^{2r} \|\hat{u}(\omega)\|_{s,\omega,\Gamma,*}^2 \, d\omega \right)^{1/2}.
	$$
\end{Definition}
As in the case of standard Sobolev spaces, these spaces are independent of the choice of $\alpha_i$ and $\varphi_i$ for $|s|\leq 1$. We refer to \cite{gimperlein2017priori} and \cite{Ha-Duong03a} for more details.
\\The weak formulation of \eqref{eq:hypersingular_equation} is based upon the bilinear form
\begin{align}\label{eq:Bilin}
 w(\psi, \Psi) := \int_{\mathbb{R}^+ \times \Gamma} (\mathcal{W} \psi)(t,x) \, \partial_t \Psi(t,x) \, d_\sigma t \, ds_x,
\end{align}
involving the weighted measure \(d_\sigma t := e^{-2\sigma t} dt\). 
\begin{theorem}[\cite{gimperlein2017priori,Ha-Duong03a}]\label{mapthm}
	Let \(r \in \mathbb{R}\). The hypersingular operator \(\mathcal{W}\) and its inverse $\mathcal{W}^{-1}$ are continuous mappings:
	\begin{align}\label{eq:W_Mapping}
			\mathcal{W} :& H_\sigma^{r+1}(\mathbb{R}^+,\widetilde{H}^{\frac{1}{2}}(\Gamma)) \to H_\sigma^{r}(\mathbb{R}^+,H^{-\frac{1}{2}}(\Gamma)),\\
		\mathcal{W}^{-1}:& H_\sigma^{r+1}(\mathbb{R}^+,\widetilde{H}^{-\frac{1}{2}}(\Gamma)) \to H_\sigma^{r}(\mathbb{R}^+, H^{\frac{1}{2}}(\Gamma))\nonumber.
	\end{align}
	Moreover, the bilinear form $w(\cdot, \cdot)$ is weakly coercive, i.e.,
	$$
	C(\sigma) \|\psi\|_{0,\frac{1}{2},\Gamma,*}^2 \leq w(\psi,\psi).
	$$
\end{theorem}
Note that the domains and ranges of $\mathcal{W}$ and $\mathcal{W}^{-1}$ do not coincide.
A variational formulation of \eqref{eq:BI1} can be written as:\\
Find $\psi \in V:=H^1_\sigma(\mathbb{R}^+,\widetilde{H}^{1/2}(\Gamma))$ such that for all $ \Psi \in V$:
\begin{align}\label{eq:varW}
	 w(\psi,\Psi) = \int_{\mathbb{R}^+ \times \Gamma} f\partial_t\Psi \,d_\sigma t \, ds_x .
\end{align}
The weak formulation \eqref{eq:varW} is well-posed in $V$ \cite{gimperlein2017priori,gimperlein2018time}.  
Due to the mapping properties of $\mathcal{W}$, as shown in Theorem \ref{mapthm}, we obtain the following corollary.
\begin{corollary}\label{col:Stabi}
Let \( f \in H^2_\sigma(\mathbb{R}^+, H^{-1/2}(\Gamma)) \) and let \( \psi \in V \) denote the unique solution of~\eqref{eq:varW}. Then, the following stability estimate hold:
	$$
	\|\psi\|_{1,\frac{1}{2}, \Gamma} \leq C_{stab}(\sigma,\omega) \|f\|_{2, -\frac{1}{2}, \Gamma}.$$
\end{corollary}
\subsection{Galerkin Discretization}
The corresponding Galerkin discretization of \eqref{eq:varW} in a finite-dimensional subspace $V_{\Delta t,h}~\subset~V$ reads:
\\Find $\psi_{\Delta t,h} \in V_{\Delta t,h}$ such that for all $ \Psi_{\Delta t,h} \in V_{\Delta t,h}$:
\begin{align}\label{eq:var2W}
	w(\psi_{\Delta t,h},\Psi_{\Delta t,h}) =  \int_{\mathbb{R}^+ \times \Gamma} f \Psi_{\Delta t,h} \,d_\sigma t \, ds_x .  
\end{align}
Then, \eqref{eq:var2W} is also well-posed in $V_{\Delta t,h}$, and we inherit the stability result from Corollary~\ref{col:Stabi}, leading to the following estimate. 
\begin{corollary}
Let \( f \in H^2_\sigma(\mathbb{R}^+, H^{-1/2}(\Gamma)) \) and let \( \psi_{\Delta t,h}~\in~V_{\Delta t,h} \) denote the unique solution of~\eqref{eq:var2W}. Then, it holds:
	$$
	\|\psi_{\Delta t,h}\|_{1,\frac{1}{2}, \Gamma} \leq C_{stab}(\sigma,\omega) \|f\|_{2, -\frac{1}{2}, \Gamma}.
	$$
\end{corollary}
The constant $C$ in the above corollary depends solely on $\sigma$ and $\omega$, and not on $h$ or $\Delta t$.\\ 
Correspondingly, one obtains convergence and best approximation of the proposed method; see, e.g., \cite{bamberger86,gimperlein2018time}. 
\begin{remark}\label{rem:SLP}
Alternative formulations based upon a single-layer potential ansatz and the adjoint double-layer operator \(\mathcal{K}'\) 
lead to a second-kind boundary integral equation:
{\((-\tfrac{1}{2}Id + \mathcal{K}')\varphi = f\)}. 
Even though the continuous variational formulation is stable, as it is equivalent to the formulation in terms of $\mathcal{W}$, 
the stability of its Galerkin approximation, see \eqref{eq:var2K} below, is not guaranteed (see Appendix §~\ref{sec:Kprime}). 

\end{remark}
\begin{remark}\label{rem:Colo}
In practice, collocation methods based upon the 
representation formula are frequently used due to their simplicity \cite{stutz2008stabilitatsverhalten}. They were extensively studied in the 1990s, for example by P.J.~Davies and D.B.~Duncan \cite{davies2004stability}. 
However, collocation approaches for wave propagation frequently exhibit instabilities, which require careful consideration especially in demanding simulation scenarios. {In spite of empirical strategies to improve the stability, no reliable solution is currently available.} 
In this way, unlike in the frequency-domain, time-dependent boundary element methods based upon collocation have not been widely used. 
\end{remark}

\section{Discretization and Algorithmic Considerations} \label{sec:Algo}
\subsection{Discretization}
We approximate $\Gamma$ by a quasi-uniform mesh {\(\mathcal{T}_S = \bigcup_{i=1}^{N_\Gamma} \Gamma_i\)}, consisting of triangular elements~$\Gamma_i$.  
The time-domain \(\mathbb{R}^+\) is uniformly partitioned into subintervals \(I_n = (t_{n-1}, t_n]\) of size \(\Delta t\), giving rise to a temporal mesh {$\mathcal{T}_T= \bigcup_{n=1}^{N_t} I_n$}.  
On \(\mathcal{T}_S\), we consider polynomial basis functions~\(\{\varphi_i^p\}\) of degree~\(p\), and on \(\mathcal{T}_T\), polynomials~\(\{\beta_n^q\}\) of degree~\(q\). These form a basis for the discrete spaces:
\begin{align*}
	V_h^p &= \{ \phi\in C(\Gamma;\mathbb{R}) :\,\phi|_{\bar{\Gamma}_i} \in \mathbb{P}^p(\Gamma_i),\ \text{continuous and } \phi|_{\partial \Gamma}=0\ \text{if } p\geq 1 \}, \\
	V_{\Delta t}^q &= \{ \Phi\in C(\mathbb{R}^+;\mathbb{R}) :\,\Phi|_{\bar{I}_n} \in \mathbb{P}^q(I_n),\ \text{continuous and }\Phi(0)=0\ \text{if } q\geq 1 \},
\end{align*}
where \(\mathbb{P}^p(\Gamma_i)\) and \(\mathbb{P}^q(I_n)\) are the spaces of polynomials of degree at most \(p\) and \(q\) on the spatial and temporal elements \(\Gamma_i\) and \(I_n\), respectively.
A discrete approximation space is defined as the tensor product $V_{\Delta t,h}^{p,q} = V_h^p \otimes V_{\Delta t}^q$, 
associated with the space-time mesh $\mathcal{T}_{S,T} = \mathcal{T}_T \times \mathcal{T}_S = \bigcup_{n,i} \square_{n,i}$, where each space-time element is of the form \(\square_{n,i} = I_n \times \Gamma_i\). A basis of $V_{\Delta t,h}^{p,q}$ is given by the tensor products of the basis functions $\beta_n^q(t)$ and $\varphi_i^p(x)$ in time and space.
\subsection{Implementation of the hypersingular operator $\mathcal{W}$}
We briefly review the discretization of $\mathcal{W}$. It will be convenient to introduce the retarded time \( t' = t - |x - y| \). As introduced by \cite{terrasse93} and investigated by \cite{joly2017mathematical} we set for computations \( \sigma = 0 \).\\
Using the subspace $V^{1,1}_{\Delta t,h} \subset V$, the ansatz functions take the form:
\begin{equation}\label{eq:AnsatzF}
	\psi_{\Delta t,h}(t,x) = \sum_{m=1}^{N_t} \sum_{i=1}^{N_s} c_{i}^{m} \beta_m^{1}(t)\varphi_i^1(x),
\end{equation}
{where $N_s$ denotes the number of spatial degrees of freedom, which for piecewise linear
basis functions coincide with the vertices of the triangular surface elements $\Gamma_i\subset \mathcal{T}_S$.} Furthermore, $\beta^1_m(t)$ denotes the hat function on $\mathcal{T}_T$ defined by
$
\beta^1_m(t)= (\Delta t)^{-1} \left( (t-t_{m-1}) \gamma^{m}(t) - (t-t_{m+1}) \gamma^{m+1}(t) \right),
$
satisfying $\beta^1_m(t_m)=1$ and $\beta^1_m(t_j)=0$ for all $j\neq m$, where
$\gamma^j(t)$ is the characteristic function of the interval $I_j$. For the time derivative of the test functions we use: 
\begin{align}\label{eq:TestF}
	\dot\Psi_{\Delta t,h}(t,x) = \gamma^n(t) \varphi_j^1(x), \quad  n = 1, \dots, N_t,\quad  j = 1, \dots, N_s.
\end{align}
Based upon the representation of the hypersingular operator \( \mathcal{W} \) from \eqref{eq:hypersingular_equation}, we obtain for the left-hand side of \eqref{eq:var2W}:
\begin{align}\label{eq:Whad}
	\int_{\mathbb{R}^+\times \Gamma} (\mathcal{W} \psi)\ \partial_t \Psi \ dt \, ds_x\  = &\frac{1}{2 \pi} \int_{0}^{\infty} \int_{\Gamma \times \Gamma}
	\Big\lbrace \frac{-n_x \cdot n_y}{|x-y|} \dot \psi(t',y) \ddot \Psi(t,x)
	\\ & + \frac{(\operatorname{curl}_{|\Gamma} \psi )(t',y) \cdot ( \operatorname{curl}_{|\Gamma} \dot\Psi)(t,x)}{|x-y|} \Big\rbrace ds_y\ ds_x \ dt.\nonumber
\end{align}
Using $\mathcal{N}=N_tN_s$ ansatz and test functions defined by \eqref{eq:AnsatzF}, \eqref{eq:TestF}, the expression \eqref{eq:Whad} for the space-time Galerkin matrix leads to the following matrix formulation of the Galerkin discretization \eqref{eq:var2W}:
\begin{equation}\label{eq:Discr}  
	\renewcommand{\arraystretch}{0.5} 
	\underbrace{
		\begin{pmatrix}
			W^0 & 0 & 0 & 0 & \cdots \\
			W^1 & W^0 & 0 & 0 & \\
			W^2 & W^1 & W^0 & 0 & \\
			W^3 & W^2 & W^1 & W^0 & \cdots \\
			\vdots &  &  & \vdots & \ddots \\
	\end{pmatrix}}_{=:\bm{\mathcal{W}} \in \mathbb{R}^{\mathcal{N} \times \mathcal{N}}}
	\underbrace{
		\begin{pmatrix}
			c^1 \\
			c^2 \\
			\vdots \\
			\vdots \\
			c^{N_t}
	\end{pmatrix}}_{=:\mathcal{C} \in \mathbb{R}^\mathcal{N}}
	=
	\underbrace{
		\begin{pmatrix}
			f^1 \\
			f^2 \\
			\vdots \\
			\vdots \\
			f^{N_t}
	\end{pmatrix}}_{=:\mathcal{F} \in \mathbb{R}^\mathcal{N}}.
\end{equation}
As detailed in Equations \eqref{eq:A}, \eqref{eq:B} in the Appendix, the assembly of each block $W^n$ includes double integrals over three light cones \( E_{n-m} \), \( E_{n-m-1} \), and \( E_{n-m-2} \). Inner integrals in \eqref{eq:A} and \eqref{eq:B} are evaluated using a composite hp-quadrature, while outer integrals are computed via standard Gauss quadrature; see \cite{gimperlein2016adaptive} and \cite{ostermann2010numerical}. For bounded surfaces $\Gamma$, the light-cone contributions $E_\ell$ vanish for 
$\ell~>~\lceil~\mathrm{diam}~(\Gamma)~/\Delta t~\rceil$. Hence, only a finite number of matrices $W^n$ needs to be assembled.
The system matrix $\bm{\mathcal{W}}$ is a lower block Toeplitz matrix, reflecting the causality of the sound radiation problem and the global tensor product space-time mesh. The coefficient vectors \( c^m = (c_1^m, \ldots, c_{N_s}^m) \) represent the discrete solution at the \( m \)-th time step. The right-hand side vectors
\[
f^m = \frac{\Delta t}{2} \boldsymbol{I} \left({f}^{m-1} + {f}^m \right)
\]
are obtained from~\eqref{eq:var2W}, where \(\boldsymbol{I}\) is the spatial mass matrix for linear basis functions. Further details regarding the assembly of the system matrix \( \bm{\mathcal{W}} \) and the derivation of the right-hand side block vector \( \mathcal{F} \) are provided in the Appendix and in \cite{gimperlein2018time}.
\subsection{Marching-on-in-time (MOT) scheme}
The lower block  triangular structure of the space-time system matrix \(\bm{\mathcal{W}}\) in \eqref{eq:Discr}, combined with the invertibility of the diagonal blocks \(W^0\), leads to an efficient explicit solution strategy using a block forward substitution. This yields the classical marching-on-in-time (MOT) shown in Algorithm \ref{Alg:MOT}, where the solution coefficients \(c^n\) at time step \(n\) are recursively computed by \eqref{eq:MOT}.
{\scriptsize
	\renewcommand{\baselinestretch}{0.9}
\begin{algorithm}
	\caption{Marching-on-in-time algorithm}\label{Alg:MOT}
	\begin{algorithmic}[1]
		\For{$n = 1$ to $N_t$}
		\If{$n \geq \left\lceil \mathrm{diam}(\Gamma)/\Delta t \right\rceil +2$}
		\State Set $W^{n-1} \gets 0$
		\Else
		\State Compute and store $W^{n-1}$ (see \eqref{eq:A} and \eqref{eq:B})
		\EndIf
		\State Compute the RHS-Block: $f^n =\frac{\Delta t}{2} \boldsymbol{I} \left({f}^{n-1} + {f}^n \right)$
		\State Solve the system {\setlength{\abovedisplayskip}{0pt}   
			\setlength{\belowdisplayskip}{0pt}   
			\begin{align}\label{eq:MOT}
				W^0 c^n =  f^n 
				- \sum_{m=1}^{n-1} W^{n-m} c^m
			\end{align}
		}
		\State Store $c^n$
		\EndFor
	\end{algorithmic}
\end{algorithm}
}\\
The linear system \eqref{eq:MOT} is solved using the preconditioned conjugate gradient (PCG) method with an incomplete Cholesky factorization as a preconditioner.
\begin{remark}
	For details on assembling the adjoint double-layer operator \(\mathcal{K}'\) and its right-hand side, see~\cite{banz2016time}. The system matrix has a structure similar to the \(\mathcal{W}\) operator, enabling an analogous MOT formulation.
\end{remark}
\section{Validation and Numerical Stability Analysis}
\label{sec:ValiStabi}
\begin{figure}[htbp]
    \centering
    \begin{subfigure}[b]{0.49\textwidth}
        \centering
        \includegraphics[width=0.95\linewidth]{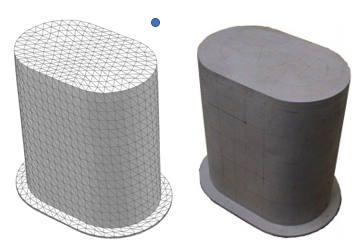}
        \caption{Oval Principle Gearbox Housing }
        \label{fig:opg_discr}
    \end{subfigure}
    \hfill
    \begin{subfigure}[b]{0.49\textwidth}
        \centering
        \includegraphics[width=0.95\linewidth]{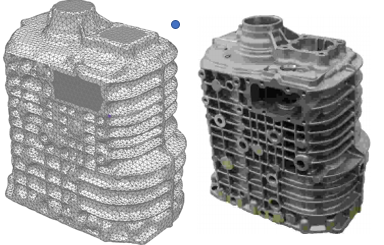}
        \caption{ZF Vehicle Gearbox Housing }
        \label{fig:zfg_discr}
    \end{subfigure}
 
    \caption{Test objects used for validation and stability analysis.}
    \label{fig:test_objects}
\end{figure}

In this section, we numerically investigate the proposed time-domain boundary element formulation concerning its convergence and its stability 
for long-time simulations. \\
Building upon our previous work (e.g., \cite{gimperlein2018time} and \cite{schneider2025stability}), we extend the analysis beyond simple test geometries to complex representative real-world gearbox housings: an oval principle housing (OPG) and a heavily ribbed ZF-gearbox housing (ZFG) from the automotive sector illustrated in Figures~\ref{fig:opg_discr} and \ref{fig:zfg_discr}.
The OPG from Honsel GmbH and Co KG, Meschede, is an 
oval-shaped housing specifically manufactured for dynamic and acoustic investigations, representing the simplest form of a gearbox housing with flat and curved surfaces and a flange. The simple design of the OPG keeps the structural complexity at a manageable level, making it suitable for fundamental studies such as convergence investigations. 
The ZFG is the ZF 6S850 Ecolite gearbox housing from ZF Friedrichshafen AG, designed for light commercial vehicles; its box-shaped structure tends to radiate sound more strongly \cite{graf2007validierung}. It represents a worst-case scenario for numerical investigations due to its highly ribbed geometry.
Throughout all numerical experiments in this section, the speed of sound is normalized to $c=1$, since it only serves as a temporal scaling parameter.
\subsection{Validation}\label{sec:Vali}
\begin{figure}
	\centering
\includegraphics[scale=0.7]{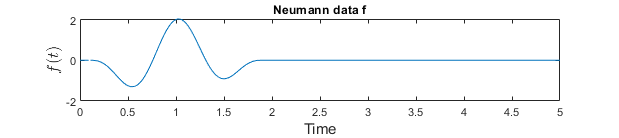}
	\caption{Neumann data $f$ for validation and stability analysis.}
	\label{fig:neumannexcitationovergamma}
\end{figure}

\begin{figure}[htb]
	\centering
	\begin{minipage}[t]{0.28\textwidth}
		\centering
		\includegraphics[width=\linewidth]{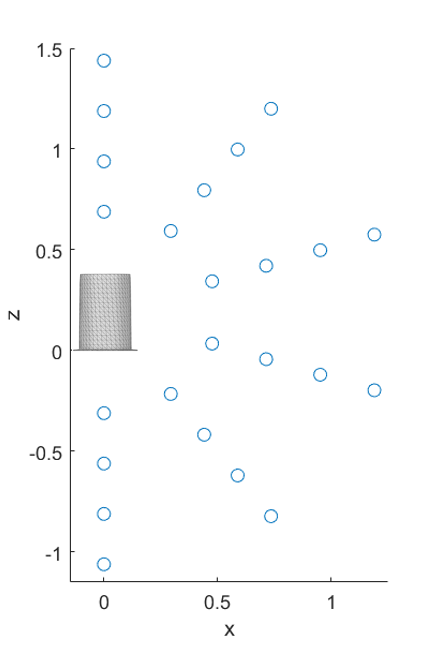}
		\caption{Mic. points.}
		\label{fig:micnodes}
		\vspace{0.02\textwidth}%
	\end{minipage}\hspace{0.02\textwidth}%
	\begin{minipage}[t]{0.55\textwidth}
		\centering
		\includegraphics[width=\linewidth]{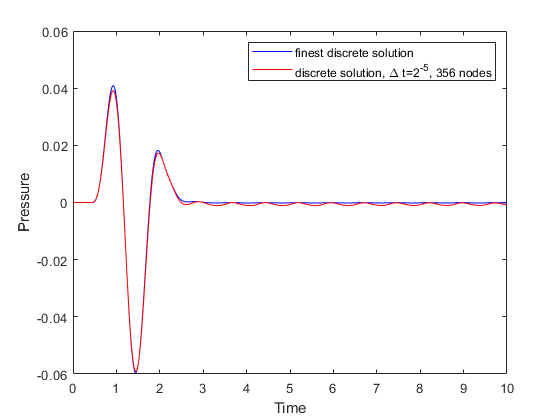}
		\caption{Pressure at \( x_1=(0,0,-0.32) \).}
		\label{fig:pressure}
	\end{minipage}
\end{figure}

\begin{figure}
	\centering
\includegraphics[scale=0.6]{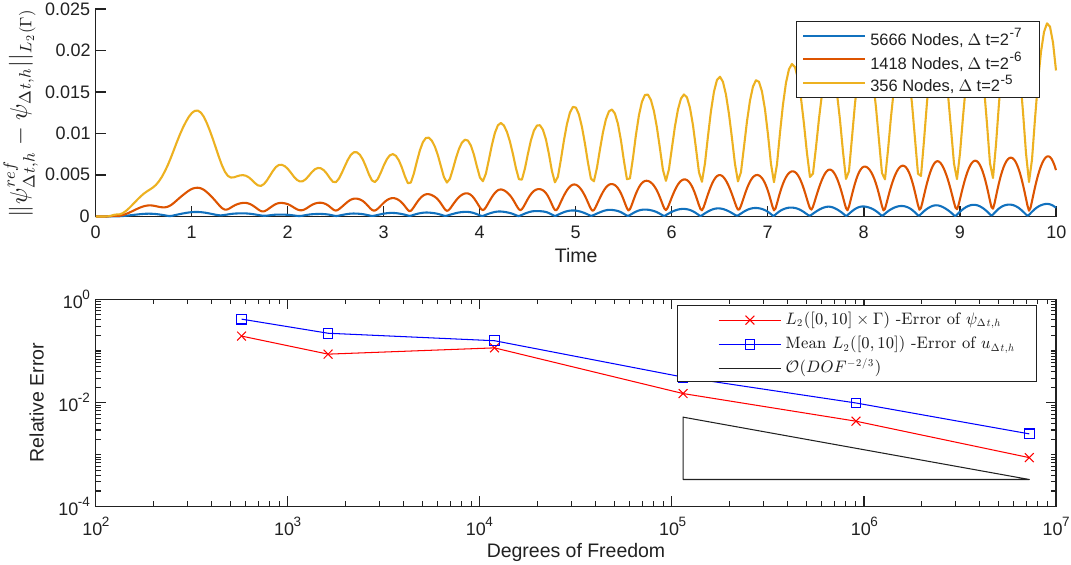}
	\caption{Relative $L_2$ error of density and sound pressure.}
	\label{fig:Conv}
\end{figure}

To assess the accuracy of the numerical determined density $\psi_{\Delta,h}$ as solution of the Galerkin approximation \eqref{eq:var2W} and the resulting discrete sound pressure
$u_{\Delta t,h}$ defined by \eqref{eq:Repr}, we examine the numerical results on the OPG from Figure~\ref{fig:opg_discr}.\\
Extending the study in \cite{schneider2025stability}, which considered the simple geometry of a unit sphere, we choose the Neumann data as follows:
\begin{align}\label{eq:rhs}
f(t) =&\Bigl[
\frac{t}{2}\left(1 + \cos\!\left(\tfrac{\pi (1-t)}{R}\right)\right)
- \frac{\pi (1-t)}{R}
   \sin\!\left(\tfrac{\pi (1-t)}{R}\right)\Bigr]\nonumber \\
   &\left(1 + \cos\!\left(\tfrac{\pi (1-t)}{R}\right)\right)H(R-|1-t|),
\end{align} 
with  $R=0.9$ for times up to $T=10$. There, $H(\cdot)$ denotes the Heaviside function, which restricts the boundary excitation to the time interval $t \in [0.1,\, 1.9]$. Outside this interval, the Neumann condition vanishes, as illustrated in Figure~\ref{fig:neumannexcitationovergamma} for the first 5 seconds.
\\
Seven space–time meshes of the OPG were considered. The spatial meshes of the test object are refined from $24$ up to $45{,}132$ elements ($14$ to $22{,}658$ nodes), while $\Delta t$ is decreased from $2^{-2}$ to {$2^{-8}$}. The resulting time-space ratio is approximately $\Delta t / h \approx 0.4$. Overall, the discretizations range from {$560$ up to $5.8 \times 10^7$ }degrees of freedom {for the experiments in this subsection}.\\
In Figure~\ref{fig:pressure}, we depict the sound pressure results at $x_1=(0,0,-0.32)$, using a spatial mesh with $356$ nodes and a time step size of $\Delta t=2^{-5}$.
The microphone points $x_i,\,i=1,...,24$ are positioned in $\mathbb{R}^3\backslash \overline{\Omega}$ on a semicircle around the OPG at a distance of 0.25 in space (see Figure~\ref{fig:micnodes}). \\
To analyze the convergence of the computed densities and sound pressures, we evaluate the error in the space–time norm:
\begin{align*}
\frac{||\psi^{ref}_{\Delta t,h} - \psi_{\Delta t,h}||_{L_2([0, 10]; L_2(\Gamma))}}{||\psi^{ref}_{\Delta t,h}||_{L_2([0, 10]; L_2(\Gamma))}} 
\quad \text{and} \quad 
\frac{||u^{ref}_{\Delta t,h}(t, x_i) - u_{\Delta t, h}(t, x_i)||_{L_2([0, 10])}}{||u^{ref}_{\Delta t,h}(t, x_i)||_{L_2([0, 10])}},
\end{align*}
 relative to the solution $\psi^{ref}_{\Delta t,h}$ and $u^{ref}_{\Delta t,h}(t, x_i)$ computed on the finest mesh.
The overall pressure error is defined as the mean of the individual errors across all 24 microphone positions. Figure~\ref{fig:Conv} illustrates the relative errors of $u_{\Delta t,h}$ and $\psi_{\Delta t,h}$ with respect to the number of degrees of freedom. The comparatively larger errors observed for the three coarsest meshes can be attributed to geometric approximation effects, since in these discretizations the flange of the OPG is not considered. Therefore, we restrict the estimation of convergence rates to the three finest meshes. The resulting convergence rate is approximately $-0.61$ for the sound pressure and $-0.68$ for the density. The rate for the density $\psi_{\Delta t,h}$ is in excellent agreement with the theoretically expected value of $-2/3$, as predicted by Theorem~18 in \cite{gimperlein2018boundary} and, for instance, demonstrated in~\cite{schneider2025stability} for the unit sphere, despite the presence of sharp edges in the geometry.

\subsection{Stability Analysis}\label{sec:Stabi}
In this section, we assess the stability of the discretized variational formulation \eqref{eq:var2W} of the $\mathcal{W}$ operator. 
The results are compared with those from two alternative approaches, namely
\begin{itemize}
\item a discretized problem based on a second-kind integral formulation \eqref{eq:var2K};
\item a semi-discretized collocation approach \eqref{eq:Colo}.
\end{itemize}
The stability analysis is carried out for a fixed spatial grid size $h$ as $\Delta t \to 0$. 
For constant ratios $\Delta t / h$, the stability is investigated in Section~\ref{sec:Vali} for the OPG while simultaneously decreasing both $h$ and $\Delta t$.
We consider both gearbox housings depicted in Figure~\ref{fig:test_objects}. Their surfaces are extended to closed geometric domains with boundary $\Gamma$, as illustrated in Figure~\ref{fig:zfg_discr}. For the Neumann data, we use \eqref{eq:rhs} over the time interval $[0,6800)$.
\begin{table}[ht]
\centering
\renewcommand{\arraystretch}{0.5} 
\begin{tabular}{l c c c }
\hline
Object & Elements & Nodes & $h$  \\
\hline
OPG & 2,832 & 1,418 & $2.2\times10^{-2}$ \\
ZFG & 28,698 & 14,351 & $1.0\times10^{-2}$  \\
\hline
\end{tabular}
\caption{Spatial discretizations for stability investigations.}
\label{tab:discretization}
\end{table}
The spatial discretizations are described in Table~\ref{tab:discretization}. Selecting time-space ratios $\Delta t / h$ from 3.0 halved five times to 0.1 yields $\Delta t$ values from $2^{-3}$ to $2^{-9}$.
Consequently, the total number of time steps $N_t$ varies from 61,608 to 4,690,500 (see Table~\ref{tab:Stabi_OPG} and \ref{tab:Stabi_ZFG}).
\\
We obtain the densities on $\Gamma$, evaluate the resulting sound pressure at $x_0=(0.5,0,0)$ (blue dot in Figure~\ref{fig:test_objects}), and consider a simulation as unstable if the computed density or sound pressure exceeds 100 times the maximum value observed in a stable reference simulation.\\
\textbf{OPG:}\\
The results calculated for the oval principle housing are summarized in Table~\ref{tab:Stabi_OPG}. It can be seen that the collocation method yields unstable results for small space-time ratios. Therefore, the time step size was further refined between $\Delta t/h = 0.19$ and $\Delta t/h = 0.10$, and instabilities were observed in all simulations using ratios $\Delta t/h < 0.19$. Similar behavior has been documented in previous studies based upon the example of a unit sphere, see~\cite{schneider2025stability}. \\
The discretized formulation of the adjoint double-layer operator $\mathcal{K}'$ yields stable results for all simulations. This is particularly remarkable when compared to results on the unit sphere~\cite{schneider2025stability}, considering the increased geometric complexity introduced by sharp edges.\\
The same holds true for the hypersingular operator, for which stability is generally expected regardless of the $\Delta t/h$ ratio.
\begin{table}[ht]
	\centering
			\renewcommand{\arraystretch}{0.5} 
	\begin{tabular}{ c|l c r|  |p{1cm}|p{1cm}|p{1cm}| }
		No.	  &  ${\Delta t / h}$& $\Delta t$ & \# timesteps  & $\mathcal{W}$ & $\mathcal{K}'$ & Colloc.   \\ 
		 \hline  \addlinespace[1pt]
		1 & 3     & $2^{-3.18}$ & 61,608    & \cellcolor{green!25}\checkmark & \cellcolor{green!25}\checkmark &  \cellcolor{green!25}\checkmark \\  
		2 &1.5   & $2^{-4.18}$ & 123,216    &\cellcolor{green!25}\checkmark  & \cellcolor{green!25}\checkmark &  \cellcolor{green!25}\checkmark \\  
		3 &0.75   & $2^{-5.18}$ & 246,500    &\cellcolor{green!25}\checkmark & \cellcolor{green!25}\checkmark &  \cellcolor{green!25}\checkmark  \\  
		4 & 0.38   & $2^{-6.17}$ &  489,668
		&\cellcolor{green!25}\checkmark & \cellcolor{green!25}\checkmark &  \cellcolor{green!25}\checkmark  \\  
		5 & 0.19 & $2^{-7.17}$ & 979,268
		&\cellcolor{green!25}\checkmark & \cellcolor{green!25}\checkmark &  \cellcolor{green!25}\checkmark  \\ 
		{6} & 0.16 & $2^{-7.41}$ & 1,156,068
		&\cellcolor{green!25}\checkmark & \cellcolor{green!25}\checkmark &  \cellcolor{red!25}\ding{55} \\ 
				{7} & 0.13 & $2^{-7.71}$ & 1,423,784
		&\cellcolor{green!25}\checkmark & \cellcolor{green!25}\checkmark &  \cellcolor{red!25}\ding{55}  \\ 
				{8} & 0.11 & $2^{-7.95}$ &1,681,436
		&\cellcolor{green!25}\checkmark & \cellcolor{green!25}\checkmark &  \cellcolor{red!25}\ding{55}  \\ 
				{9} &  0.10 & $2^{-8.09}$ & 1,852,796
		& \cellcolor{green!25}\checkmark & \cellcolor{green!25}\checkmark &  \cellcolor{red!25}\ding{55}  \\  
	\end{tabular}
	\caption{Stability Results for Different Time-Grid Resolutions (OPG)}
	\label{tab:Stabi_OPG}
\end{table}\\
\textbf{ZFG:}\\
To extend these investigations to a more complex structure, we consider the ZFG.
We did not achieve stable results with the collocation method for any of the chosen time-space ratios as can be seen in Table~\ref{tab:Stabi_ZFG}. The same applies to the adjoint double-layer potential. Notably, simulations using these two methods developed large oscillations shortly after the start of the computation.
In contrast, the first-kind integral equation ($\mathcal{W}$-operator) yielded stable results in all simulations, even for very large numbers of time steps.
\begin{table}[ht]
	\centering
				\renewcommand{\arraystretch}{0.5} 
	\begin{tabular}{ c|l c r|  |p{1cm}|p{1cm}|p{1cm}| }
	No.	  &  ${\Delta t / h}$& $\Delta t$ & \# timesteps  & $\mathcal{W}$ & $\mathcal{K}'$ & Colloc.  \\
	\hline \addlinespace[1pt]
		1 & 3    & $2^{-4.52}$ & 156,000    & \cellcolor{green!25}\checkmark &\cellcolor{red!25}\ding{55}  &  \cellcolor{red!25}\ding{55} \\  
		2 & 1.5   & $2^{-5.52}$ & 312,000    & \cellcolor{green!25}\checkmark  &\cellcolor{red!25}\ding{55} &  \cellcolor{red!25}\ding{55} \\  
		3 & 0.75   & $2^{-6.52}$ & 624,100     & \cellcolor{green!25}\checkmark  &\cellcolor{red!25}\ding{55} & \cellcolor{red!25}\ding{55} \\  
		4 &  0.38  & $2^{-7.50}$ & 1,230,900
		&\cellcolor{green!25}\checkmark &\cellcolor{red!25}\ding{55} &  \cellcolor{red!25}\ding{55} \\  
		5 &  0.19 & $2^{-9.50}$ & 2,461,900
		&\cellcolor{green!25}\checkmark  &  \cellcolor{red!25}\ding{55}  & \cellcolor{red!25}\ding{55} \\ 
		6 &  0.10 & $2^{-9.43}$ & 4,690,500
		&\cellcolor{green!25}\checkmark &  \cellcolor{red!25}\ding{55} &  \cellcolor{red!25}\ding{55}  \\  
	\end{tabular}
	\caption{Stability Results for Different Time-Grid Resolutions (ZFG)}
	\label{tab:Stabi_ZFG}
\end{table}
\\
In summary, the collocation method leads to stable results for simple geometries and sufficiently large time-space ratios but quickly reaches its limits for more complex structures. The ansatz using a single-layer potential exhibits similar behavior; however, no clear pattern for a suitable discretization emerges, as it appears to strongly depend on the mesh on $\Gamma$ and not only on the geometric complexity.  
The numerical results reflect the absence of rigorous stability results for both discrete representations, {or conceivably also a larger influence of quadrature errors in the matrix assembly}.\\
The formulation based upon the $\mathcal{W}$-operator consistently gives stable results in all experiments, including long-term simulations with over 4.5 million time steps. This is in accordance with its proven well-posedness and stability.
\section{Experimental Setup and Numerical NVH Investigations}
\label{sec:Exp}
\begin{figure}[htbp]
        \centering
        \includegraphics[width=1\linewidth]{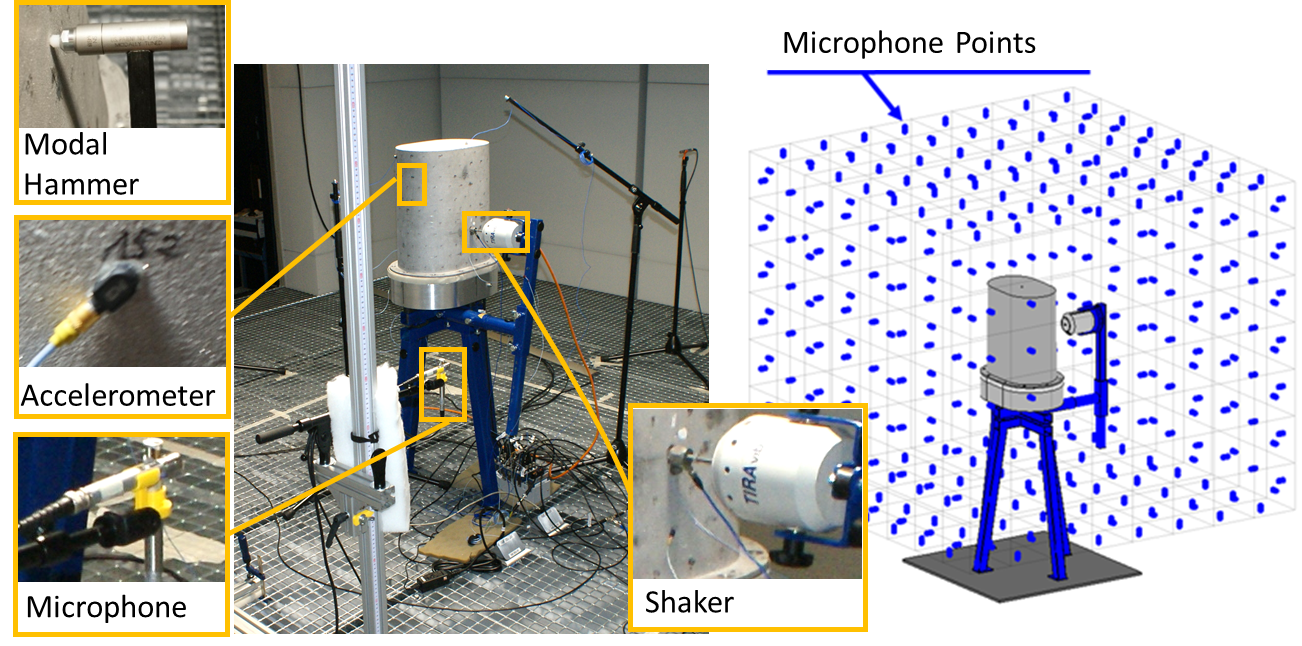}
        \caption{Test bench in anechoic chamber}
        \label{fig:TestBench}
\end{figure}

\begin{figure}[ht]
	\centering
			\renewcommand{\arraystretch}{0.5} 
	\begin{tikzpicture}[
		node distance=1.cm and 1.cm,
		box/.style={draw=black, fill=cyan!30, minimum width=2.6cm, minimum height=0.71cm, align=center},
		smallbox/.style={draw=black, fill=gray!30, minimum width=2.5cm, minimum height=0.8cm, align=center},
		every node/.style={font=\sffamily},
		arrow/.style={-{Latex}, thick},
		]
		
		\node[box] (Str) {Structural Analysis};
		\node[box, right=of Str] (Dyn) {Dynamic Analysis};
		\node[box, right=of Dyn] (Acoustic) {Acoustic Analysis};
		\node[smallbox, below=of Acoustic] (air) {\shortstack{Sound pressure\\ $u$ in $\mathbb{R}^+\!\times\!\mathbb{R}^3\backslash \Gamma$}};

		\draw[arrow] (Acoustic) -- (air);
		
		\node[above=0.3cm of Str] (femesh) {FE-Mesh of structure};
		\draw[arrow] (femesh) -- (Str);		
		\node[below=0.3cm of Str] (eig) {%
			\begin{tabular}{l}
				Eigenfrequencies \\
				Mode shapes
		\end{tabular}};
		\draw[arrow] (Str) -- (eig);
		\draw[arrow] (eig.east)  -- ++(1.54,0) -- ([xshift=-20pt]Dyn.south);

	\node[below=0.3cm of Dyn, xshift=16pt, yshift=0pt] (g) {%
	
	\begin{tabular}{l}
		Neumann data \\ $f$
		on $\mathbb{R}^+\!\times\!\Gamma$
	\end{tabular}};

		\draw[arrow] ([xshift=16pt]Dyn.south) --  ([xshift=0pt]g.north);
		
				\draw[arrow] (g.east)  -- ++(1.05,0) -- ([xshift=-20pt]Acoustic.south);
	\node[above=0.3cm of Dyn] (inputs) {
		\begin{tabular}{l}
			Excitation \\
			Structural Damping
		\end{tabular}
	};
		\draw[arrow] (inputs) -- (Dyn);
		
	\end{tikzpicture}
	\caption{NVH calculation process}
	  \label{fig:gearbox-radiation}
\end{figure}
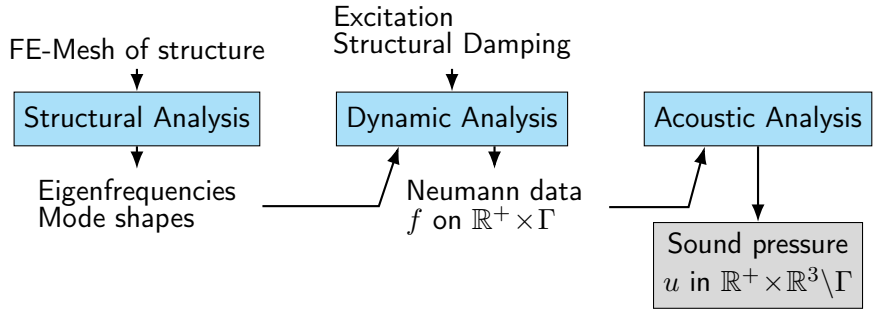

\textbf{Experimental Setup:}
Experimental validation is performed on the test objects described in Section~\ref{sec:ValiStabi}, using a test rig located in a fully anechoic chamber (Figure~\ref{fig:TestBench}).
The rig is specifically designed to ensure minimal mechanical coupling to the environment across the frequency range of interest. It comprises a sand-filled steel frame, a stiff aluminum plate whose first natural frequency lies above the highest frequency under investigation, and a foam layer, all contributing to mechanical decoupling of the test objects from the environment.\\
During both transient and steady-state tests, a complementary suite of sensors monitors the dynamic and acoustic responses. Structure-borne signals are captured by multiple miniature accelerometers and a scanning laser Doppler vibrometer, while sound radiation is measured at discrete positions using several microphones.\\
Over the course of the experiments, excitation at the test bench is provided by either an instrumented modal hammer or an electromagnetic shaker, enabling both broadband and single-frequency inputs.\\
\textbf{Numerical NVH Analysis:}
To compute the sound radiation of the various test cases examined in the anechoic chamber, the calculation workflow depicted in Figure~\ref{fig:gearbox-radiation} is used. A validated structural analysis provides modal models of the OPG and ZFG, including their mode shapes and eigenfrequencies as well as the modal dampings extracted by an experimental modal analysis. These models, together with the excitation force signal from the physical tests, are used in a reduced flexible multibody simulation to compute the surface accelerations and thus the Neumann boundary data $f$ in \eqref{eq:ori_problem_NeumannData}. The validation of the dynamic analysis, and consequently of the numerical Neumann data with respect to the physical tests, is presented in \cite{schneider2022practical}.
\section{Sound Emission of a Simplified Gearbox Housing}
\label{sec:6}
In this section, we use the $\mathcal{W}$-operator to compute the sound radiation $\psi_{\Delta t,h}$ as the solution of \eqref{eq:var2W} for the oval principal housing introduced in Section~\ref{sec:Exp}, with particular emphasis on validation against experimental measurements. 
For this purpose, the speed of sound is set to $c = 340$ m/s in the following investigations.\\
The Neumann data is obtained from a structural and dynamic simulation, as illustrated in Figure~\ref{fig:gearbox-radiation}. Two test scenarios are considered: first, a harmonic excitation where the test object is driven by a sinusoidal input; second, an impact excitation covering a broad frequency range of interest.\\
For our computations, the time step size is set to $\Delta t = 1/25{,}600$, chosen to match the sampling rate of the measurement system. 
The spatial mesh of the outer surface of the test object $\Gamma$ forms a closed hull and includes $9{,}784$ elements and $4{,}876$ nodes, resulting in a time-space relation of $c \Delta t / h \approx 1.1$. 
The lower side of the OPG flange region, which is positioned on the aluminium plate of the test bench (see Figure~\ref{fig:TestBench}), is given a zero vector as the right-hand side.
\\
The sound pressure signals $u_{\Delta t,h}(t,x_i)$ (cf.~\eqref{eq:Repr}) are obtained from $\psi_{\Delta t,h}$ at discrete points in the exterior $\mathbb{R}^3 \setminus \overline{\Omega}$, enabling direct comparison with the signals recorded by the microphones.
\subsection{Harmonic Sound Investigation}
\begin{figure}
	\centering
	\includegraphics[width=0.9\linewidth]{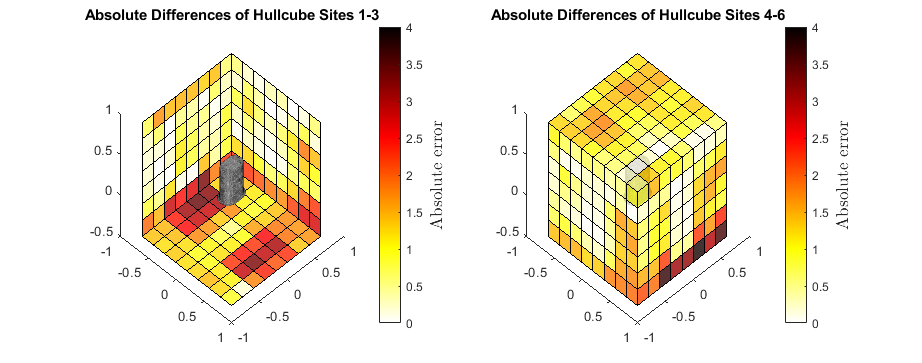}
\caption{Deviations of sound pressure levels over $\mathcal{H}$}
	\label{fig:l2errorhullcube}
\end{figure}

The test object is excited at its first natural frequency of $386$ Hz, promoting pronounced sound radiation. 
For the calculations, we consider the time range $[0,3)$ while the evaluation is limited to $[2,3)$, when a steady state has been reached.\\
Sound pressures are evaluated at $N_M=352$ uniformly distributed points $x_i$ on the cuboid surface $\mathcal{H}$ enclosing the test object (Fig.~\ref{fig:TestBench}). Each $x_i$ is the centroid of a square subarea $\mathcal{H}_i\subset \mathcal{H}$ with 20 cm edge length and $\mathcal{H}=\bigcup_i\mathcal{H}_i$.\\
Following DIN EN ISO 3744, the sound power level $L_W(u;\mathcal{H},[t_1,t_2])$ is derived from the sound pressure levels (SPL) $L_p(u_i(t);[t_1,t_2])$ at microphone positions $x_i$ by energy-based averaging over the associated subareas $\mathcal{H}_i$:
\begin{align*}
	L_p(u_i(t);[t_1,t_2]) &:= 10 \log \left[ \frac{1}{T} \frac{\| u_i(t) \|_{L_2([t_1,t_2])}^2}{u_0^2} \right],  \\
	L_W(u;\mathcal{H},[t_1,t_2])  &:= 10 \log \left[ \frac{1}{N_M} \sum_{i=1}^{N_M} 10^{0.1 L_p(u_i(t);[t_1,t_2]) }\cdot|\mathcal{H}_i| \right].
\end{align*}
with reference value $u_0=2\cdot 10^{-5}$.
\begin{figure}
	\centering
	\includegraphics[width=0.9\linewidth]{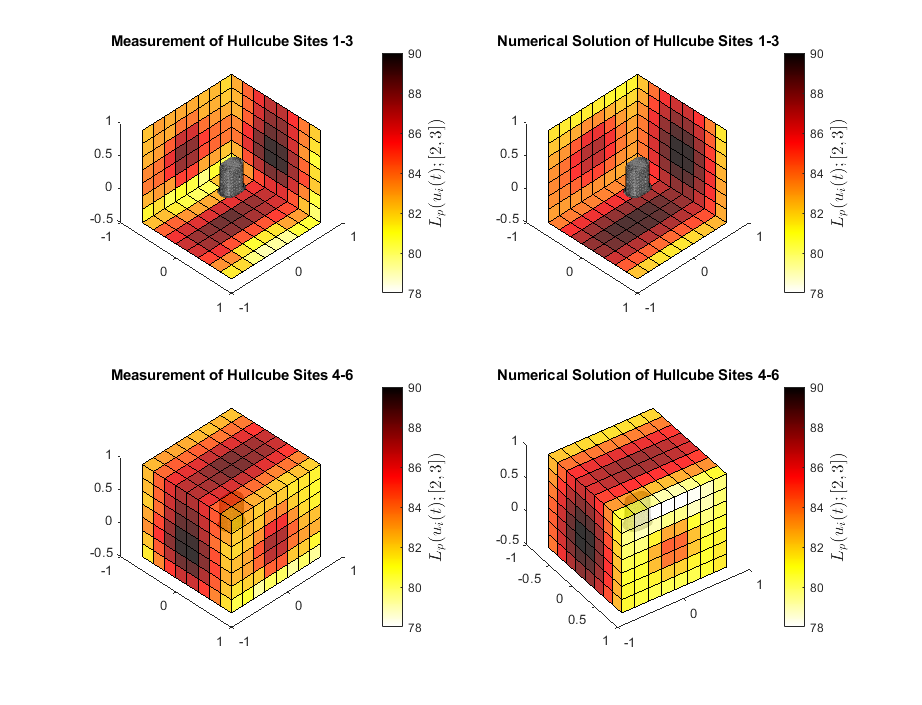}
\caption{Sound pressure levels over $\mathcal{H}$}
	\label{fig:soundpower}
\end{figure}
Figure~\ref{fig:soundpower} shows the SPL for each subarea $\mathcal{H}_i$. 
The comparison reveals a high degree of similarity between the numerical and measured sound radiation.
However, the measured pressures tend to be slightly higher on the front side of $\mathcal{H}$. Figure~\ref{fig:l2errorhullcube} therefore shows, for each surface $\mathcal{H}_i$, the SPL errors, i.e., the difference between the measured and calculated signals.
The deviations observed in the lower range of Figure \ref{fig:l2errorhullcube} can be attributed to discrepancies between the test rig configuration and the acoustic BEM mesh, as the test rig itself is not represented in the numerical simulation. We obtain $L_W(u;\mathcal{H},[2,3]) = 85.36\,\mathrm{dB}$ for the simulation and $L_W(u;\mathcal{H},[2,3])  = 85.14\,\mathrm{dB}$ for the measurement, both describing the overall sound radiation of the complete test object. The small deviation between these values demonstrates the accuracy of the numerical results in capturing the real-world acoustic behavior of the OPG.
\subsection{Transient Sound Investigation}
\begin{figure}
	\centering
	\includegraphics[width=0.8\linewidth]{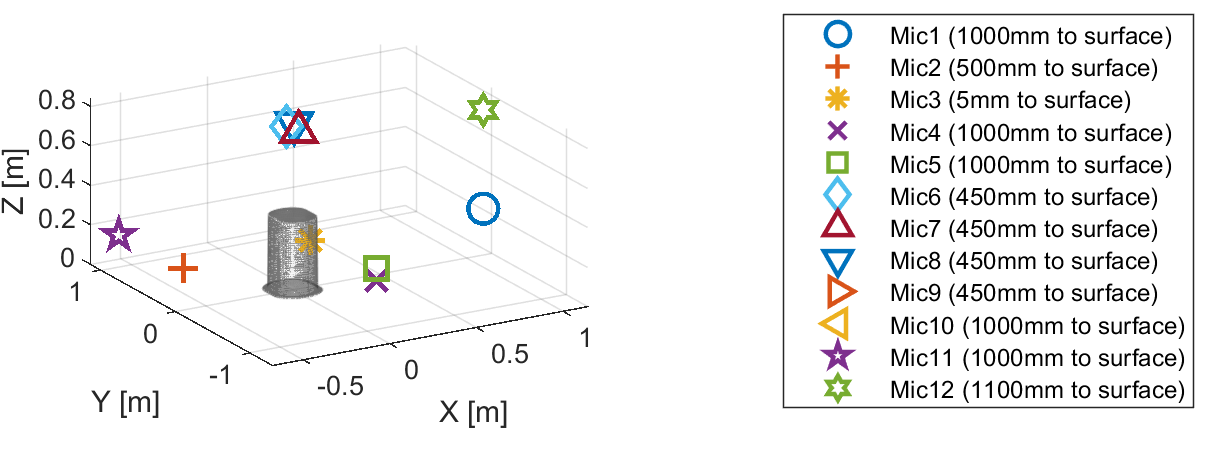}
	\caption{Microphone positions for transient OPG investigation}
	\label{fig:opgmicsfarfield}
\end{figure}

\begin{figure}
	\centering
	\includegraphics[width=0.9\linewidth]{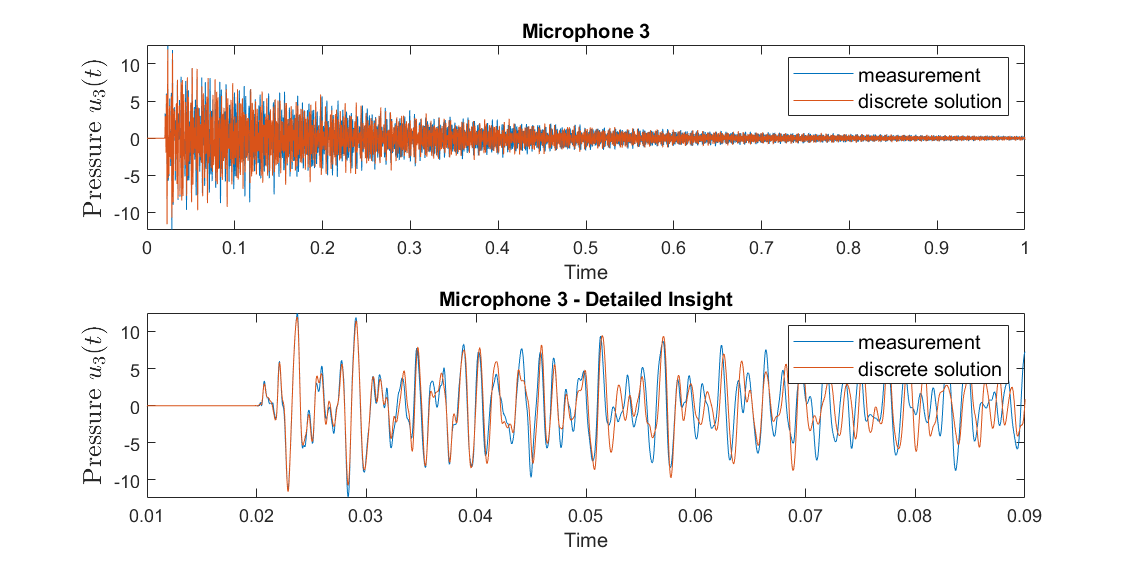}
	\caption{Sound pressure of microphone 3 in time-domain, OPG}
	\label{fig:Mic3_Farfield}
\end{figure}

\begin{figure}
	\centering
	\includegraphics[width=0.9\linewidth]{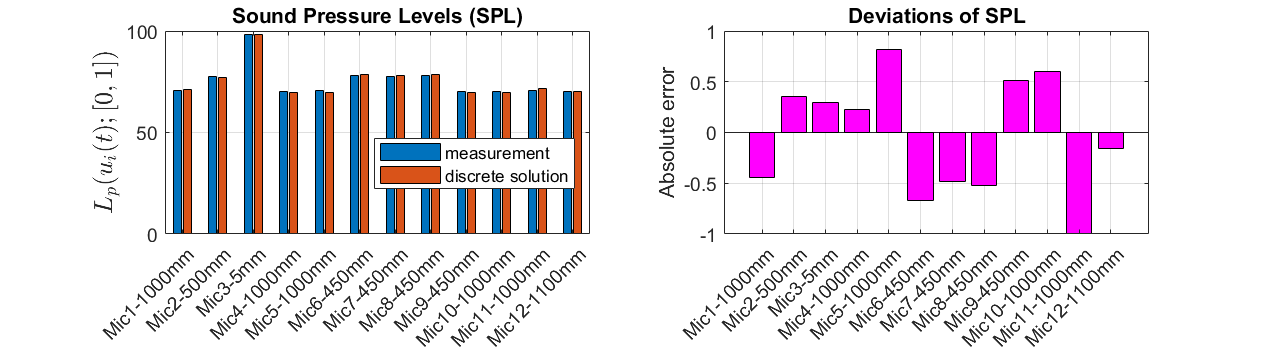}
	\caption{Deviations of sound pressure levels in time-domain, OPG}
	\label{fig:L2Error_FarField}
\end{figure}

\begin{figure}
	\centering
	\includegraphics[width=0.9\linewidth]{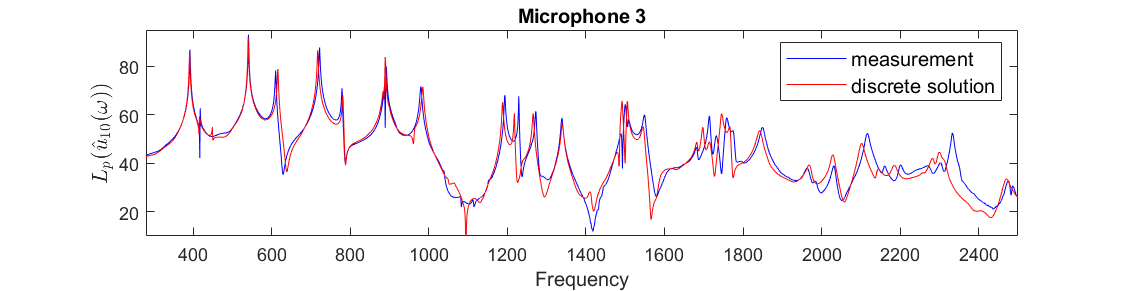}
	\caption{Sound pressure of microphone 3 in frequency-domain, OPG}
	\label{fig:Mic3_FD}
\end{figure}

In addition to the harmonic case, a transient excitation is investigated using a modal hammer with an integrated force sensor. This allows the excitation of a broad frequency range up to 2500 Hz. \\
For the calculation, the time interval is limited to $T=1$, since the structural vibrations -- and thus the Neumann boundary conditions $f$ -- have been reduced to negligible levels within this period. Sound pressures are measured and computed at $11$ points, as illustrated in Figure~\ref{fig:opgmicsfarfield}. These are located both on the surface of the test object and at positions up to 1.4~m away in the anechoic chamber.\\ Figure~\ref{fig:Mic3_Farfield} shows the temporal evolution of the experimental and numerical sound pressure signals at microphone 3, located close to the surface of the test object. 
The magnitudes of both signals agree closely during the initial time steps; however, slight phase shifts emerge after a short period. These deviations are consistent with effects observed in the dynamic analysis. Thus, the transient acoustic signals are evaluated in terms deviations between the level values. Figure \ref{fig:L2Error_FarField} presents the SPL values for each microphone position. 
Further, the deviations between measured and simulated level values are shown on the right-hand side, remaining below $1$ dB across all evaluation points $x_i$.
A frequency-domain analysis complements the time-domain results by transforming the signals into
$\hat{u}_i(\omega) := \hat{u}(\omega,x_i) = \mathcal{F}_t(u(t,x_i)), \, \omega \in \mathbb{R}$, using the Fourier transform $\mathcal{F}_t$ in time. In the numerical tests, this transformation is performed using the Fast Fourier Transform (FFT).
Figure~\ref{fig:Mic3_FD} presents the sound pressure level $\hat{u}_i(\omega)$ in the frequency range $[280,2800)$ Hz, where structural vibrations dominate the excitation. Good agreement is observed up to 2100 Hz, with amplitudes at the main peaks deviating by at most $1$ dB. 
The deviations observed up to 2500 Hz can primarily be attributed to the corresponding discrepancies in the structural and dynamic analysis (see \cite{schneider2022practical}). 
The remaining frequency signals of the 11 microphones exhibit a similarly good match. 
In conclusion, the transient signals $u(t,x_i)$ show excellent agreement between numerical and experimental results, with SPL deviations below 1\,dB in both the time and frequency-domains.
\section{Transient Sound Emission of a Complex Gearbox Housing}\label{sec:7}
\begin{figure}
	\centering
	\includegraphics[width=0.9\linewidth]{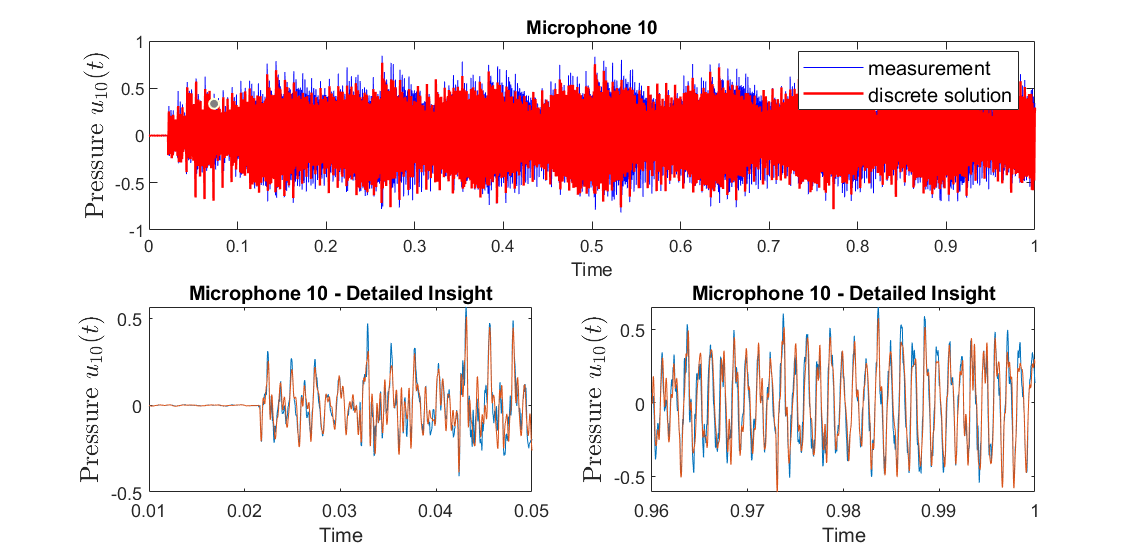}
	\caption{Sound Pressure of microphone 10 in time-domain, ZFG}
	\label{fig:mic12tc5}
\end{figure}

\begin{figure}
	\centering
	\includegraphics[width=0.9\linewidth]{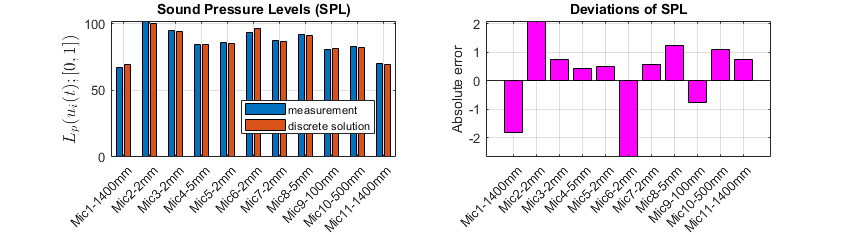}
	\caption{Deviations of sound pressure levels in time-domain, ZFG}
	\label{fig:tc5errordbtd}
\end{figure}

\begin{figure}
	\centering
	\includegraphics[width=0.9\linewidth]{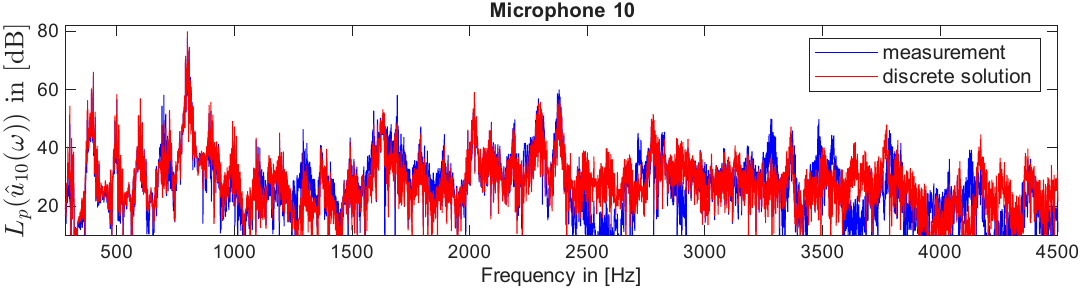}
	\caption{Deviations of sound pressure in frequency-domain, ZFG}
	\label{fig:fdsplmic111tc5}
\end{figure}

In this section, the Galerkin approximation \eqref{eq:var2W} of the first-kind integral equation \eqref{eq:BI1} is applied to compute the sound field of the complex automotive gearbox housing introduced in Section~\ref{sec:ValiStabi}. The focus is on comparing the numerical results with corresponding experimental measurements from the test bench described in Section~\ref{sec:Exp}.\\
A transient test scenario is considered, where the gearbox housing is excited by a random impact of the shaker at a specific surface location. The Neumann data is obtained by dynamic simulations illustrated in Figure~\ref{fig:gearbox-radiation}.
\begin{figure}
	\centering
	\includegraphics[width=1\linewidth]{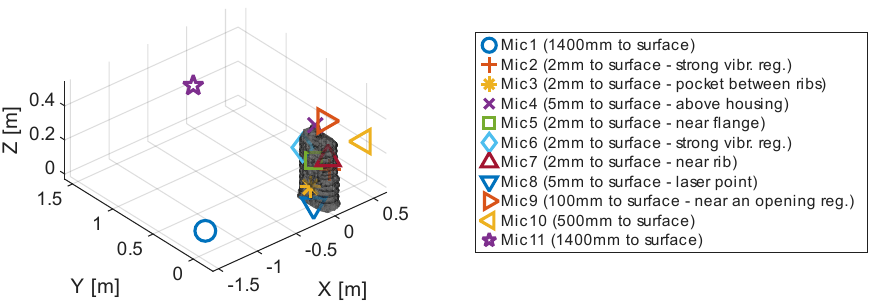}
	\caption{Microphone positions for transient ZFG investigation}
	\label{fig:zfgmics}
\end{figure}
We choose a time step size of $\Delta t = 1/51{,}200$ (according to the measurement equipment) and a mesh with 14{,}351 nodes and 28{,}698 elements, yielding $c\Delta t/h \approx 0.7$. \\
The closed regions of $\Gamma$ (Figure~\ref{fig:test_objects}) are assigned a zero right-hand side, as a foam block inside the ZFG reduces sound propagation from the interior.\\
As illustrated in Figure~\ref{fig:zfgmics}, eleven microphones are arranged around the test object, including positions close to the surface (on ribs and within pockets) as well as locations up to $1.4\,\mathrm{m}$ away.\\ Figure~\ref{fig:mic12tc5} demonstrates the agreement of the computed and measured sound pressures at microphone point~10, located 500\,mm from the test object. Remarkably, the computed and measured signals continue to agree very well even after more than 48{,}000 time steps, as shown in the detailed view at the bottom of Figure~\ref{fig:mic12tc5}, with maximum deviations remaining below $7\%$.
To evaluate the signals across all microphone positions, the SPL values are analyzed to examine the deviations. The results are presented in Figure~\ref{fig:tc5errordbtd}, showing deviations of less than $2.8$\,dB, which can be considered quite acceptable for such a complex structure.\\
In parallel to the transient investigations of the OPG in Section~\ref{sec:7}, 
the time-dependent signals are transformed into the frequency range. 
We consider the frequency interval $[300,4500)$, since this range is validated 
by the preceding structural and dynamic analyses. 
As shown in Figure~\ref{fig:fdsplmic111tc5}, the sound pressure level values 
of $\hat{u}_1(\omega)$ agree very well up to 3800\,Hz. 
In particular, in the frequency range up to 1000\,Hz, which contains the dominant 
contributions of the time signals, the magnitudes closely match. 
In this region, the maximum deviations between measurement and simulation remain 
below $0.3\,\%$ in terms of the maximum amplitudes in dB. 
The discrepancies that become visible above 2500\,Hz can be attributed to 
differences in the structural vibrations.
Moreover, the amplitudes in this higher-frequency range are approximately 
20--30\,dB lower than the dominant components, and therefore only 
marginally contribute to the overall time-domain response.
Notably, the numerical results reveal no locally increased sound pressure levels within regions such as ribs and pockets, although the Neumann data imposed strong excitation over the entire time interval up to $T=1$\,s.
\section{Conclusions}\label{sec:8}
We presented a well-posed and numerically stable space-time discretization of the hypersingular operator for the acoustic Neumann problem of the time-dependent wave equation. While previously considered only theoretically, we tested its performance for complex geometries, long trajectories with millions of time steps, and realistic Neumann boundary data. The Galerkin approximation yields a lower block-Toeplitz system, enabling efficient time-stepping solutions.\\
Long-term stability analysis confirms that the formulation remains robust even for very fine temporal resolutions and extremely long trajectories, whereas conventional {second-kind} BEM or collocation approaches can become unstable depending on the discretization or the geometry. Simulations of gearbox sound radiation show excellent agreement with anechoic chamber measurements, with harmonic deviations below $0.3$ dB and similarly accurate transient results.\\
These findings demonstrate the strong potential of the hypersingular operator based formulation for industrial applications, particularly due to its robustness with respect to geometric complexity and long time trajectories. The main limitation is the computational cost of assembling the operator, motivating future work on matrix compression and reduced-order methods to accelerate both assembly and simulation for larger problems.
\section{Appendix}
We describe briefly two additional approaches often used in practical applications: one based on the adjoint double-layer potential, and another using a collocation method with a combined layer ansatz. Since we compare our BEM with these methods, we collect all details of the discretizations we have used in the comparisons. 
\subsection{Second Kind Boundary Integral Equation}\label{sec:Kprime}
A single-layer potential ansatz $u(t,x)=S\varphi(t,x)$ 
leads to the boundary integral equation of second kind
$(-\tfrac12 \mathrm{Id}+\mathcal{K}')\varphi = \frac{\partial}{\partial n} u(t,x)$,  
where the adjoint double-layer operator is defined by
\begin{align}\label{eq:ADLP}
\mathcal{K}' \varphi := \int_{\mathbb{R}^+} \int_{\Gamma} \frac{\partial G}{\partial n_x}(t-t^*,x,y)\, \varphi(t^*,y)\, dt^*\, ds_y.
\end{align}
Following \cite{gimperlein2018time}, a variational formulation is considered in the continuous setting, for which well-posedness and stability in the corresponding anisotropic Sobolev spaces can be shown. 
The tensor-product space $V_{\Delta t,h}^{p,q}$  provides a Galerkin approximation: find $\varphi_{\Delta t,h}\in V_{\Delta t,h}^{p1,q1}\subset H^1_\sigma(\mathbb{R}^+,\widetilde{H}^{-1/2}(\Gamma))$ such that for all $\phi\in V_{\Delta t,h}^{p2,q2}\subset H^1_\sigma(\mathbb{R}^+,\widetilde{H}^{1/2}(\Gamma))$:
\begin{align}\label{eq:var2K}
\int_{\mathbb{R}^+\times \Gamma} (-\tfrac12 \mathrm{Id}+\mathcal{K}')\varphi_{\Delta t,h}\, \partial_t \phi_{\Delta t,h}\, ds_x\, d_\sigma t
=
\int_{\mathbb{R}^+\times \Gamma} f\, \partial_t \phi_{\Delta t,h}\, ds_x\, d_\sigma t.
\end{align}
For the computations, we set $\sigma = 0$ and choose $\varphi_{\Delta t,h},\, \partial_t \phi_{\Delta t,h} \in V_{\Delta t,h}^{0,0}$, which yields an MOT scheme \cite{banz2016time,ozdemir2019finite}. 

{In each time step this scheme requires the computation of a sparse Galerkin matrix, as described in \cite{banz2016time}. Like for the hypersingular operator $\mathcal{W}$, the accurate numerical evaluation of the entries is based, for the inner integral over $\Gamma$, on a geometric decomposition into light cones combined with a composite $hp$ quadrature rule. A composite Gauß quadrature is then used for the outer quadrature. As discussed in \cite{gimperlein2016adaptive,ostermann2010numerical,iccs2008} the resulting error decreases exponentially with the fourth root of the number of quadrature points.}

\subsection{Collocation-Based Approach}\label{sec:colo}
Building upon the framework used in~\cite{stutz2008stabilitatsverhalten}, we describe the solution $u \in \mathbb{R}^+ \times \mathbb{R}^3\backslash\overline{\Omega}$ by a combined layer potential ansatz $u(t,x) = -S \frac{\partial u^e}{\partial n} + D u^e$, where $u^e$ is the outer trace of $u$. Using an approximation of $u$ that is linear in time and constant in space, 
and assuming $f$ to be constant in both space and time, 
a semidiscretization leads to
\begin{align}\label{eq:Colo}
    u_h^n(x_i) \;=\; -\tfrac{1}{2}\, S_{\Delta t,h} f_h^n(x_i) 
    + D_{\Delta t,h}\, u_h^{n,+}(x_i).
\end{align}
{For the calculations presented in Section~\ref{sec:Stabi}, we use the implementation from~\cite{stutz2008stabilitatsverhalten}, where the quadrature is described. 
In particular, we note that due to the use of linear elements in this implementation, the strong singularity in \eqref{eq:Colo} vanishes \cite[Chapter~3.4]{stutz2008stabilitatsverhalten}.}

\subsection{Implementation of $\mathcal{W}$}
The implementation follows the approach in~\cite{gimperlein2018time}. From \eqref{eq:Whad}, the discrete test functions \eqref{eq:TestF}, and the trial functions \eqref{eq:AnsatzF}, we obtain
\begin{align*}
	\int_{\mathbb{R}^+\times\Gamma} \mathcal{W} \psi_{\Delta t,h}(t,x)\ \partial_t \Psi_{\Delta t,h}(t,x) \  dt\ ds_{x} 
	&=  A-B,
\end{align*}
with
\begin{align*}
		A= \sum_{m=1}^{N_t} \sum_{i=1}^{N_s} \frac{c^m_i}{2 \pi} \int_{\Gamma \times \Gamma} \frac{n_x \cdot n_y}{|x-y|} \varphi^1_{i}(y) \varphi^1_{j}(x) 
	\left( \int_{0}^{\infty} \dot \beta_1^{m}(t') \dot \gamma^{n}(t)\ dt \right) ds_y\ ds_x,
\end{align*}
and
\begin{align*}
	B =\sum_{m=1}^{N_t} \sum_{i=1}^{N_s} \frac{c_i^m}{2 \pi} &\int_{\Gamma \times \Gamma } 
	\frac{1}{|x-y|} \\
	&\int_{0}^{\infty} \beta_1^{m}(t')\ \operatorname{curl}_{|\Gamma} \varphi^1_{i}(y) \
	\gamma^{n}(t)\ \operatorname{curl}_{|\Gamma} \varphi^1_{j}(x)\ dt\ ds_y\ ds_x .
\end{align*}
With the Dirac distribution $\delta_{t_n}$ and $\dot \gamma^{n} = \delta_{t_{n-1}}-\delta_{t_n}$, we obtain for the inner time integral in $A$:
\begin{align*}
	\int_{0}^{\infty} \dot \beta_1^{m}(t') \dot \gamma^{n}(t) dt  =  - (\Delta t)^{-1} \left( \chi_{E_{n-m}}(x,y) - 2 \chi_{E_{n-m-1}}(x,y) + \chi_{E_{n-m-2}} \right) .
\end{align*}
Here, $E_l$ denotes a light cone, defined as $E_l := \lbrace (x,y) \in \Gamma \times \Gamma : t_{l} \leq |x-y| \leq t_{l+1} \rbrace \subset \Gamma \times \Gamma$ and the indicator function $\chi_{E_l}(x,y)$, which is equal $1$ if $(x,y) \in E_l$, and $0$ otherwise. Thus, we conclude
\begin{align}\label{eq:A}
	A =  \sum_{m=1}^{N_t} \sum_{i=1}^{N_s} c^{m}_i 
	\Bigg[ &- \int\limits_{E_{n-m}} \frac{(n_x \cdot n_y)(\Delta t)^{-1} 
		\varphi^1_{i}(y) \varphi^1_{j}(x)}{2 \pi |x-y|} ds_y\ ds_x \nonumber \\
	& + 2 \int\limits_{E_{n-m-1}} \frac{(n_x \cdot n_y)(\Delta t)^{-1} 
		\varphi^1_{i}(y) \varphi^1_{j}(x)}{2 \pi|x-y| } ds_y\ ds_x \nonumber \\
	& - \int\limits_{E_{n-m-2}} \frac{(n_x \cdot n_y)(\Delta t)^{-1} 
		\varphi^1_{i}(y) \varphi^1_{j}(x)}{ 2 \pi|x-y|} ds_y\ ds_x \Bigg]\ .
\end{align}
Rearranging the terms in $B$, we obtain:
\begin{align*}
	B= \sum_{m=1}^{N_t} \sum_{i=1}^{N_s} \frac{c_i^m}{2 \pi}& \int_{\Gamma \times \Gamma} 
	\frac{1}{|x-y|} \operatorname{curl}_{|\Gamma} \varphi^1_{i}(y) \operatorname{curl}_{|\Gamma} \varphi^1_{j}(x) \\
	&\int_{0}^{\infty} \beta_1^{m}(t') \gamma^{n}(t)\ dt \ ds_y\ ds_x\ .
\end{align*}
By substituting the definition of $ \beta_1^{m}$, the time integral becomes:
\begin{align*}
	&\int_{0}^{\infty} \beta_1^{m}(t') \gamma^{n}(t)\ dt\\ &\quad = (2 \Delta t )^{-1} (|x-y|^{2} - 2 |x-y| t_{n-m+1} + t_{n-m+1}^{2} ) \chi_{E_{n-m}} (x,y) \\
	&\quad\quad+ (2 \Delta t )^{-1} ( |x-y|^{2} - 2 |x-y| t_{n-m-2} + t_{n-m-2}^{2} ) \chi_{E_{n-m-2}}(x,y) \\
	&\quad\quad+ (2 \Delta t )^{-1} ((-2 |x-y|^{2} + 2 |x-y| (t_{n-m} + t_{n-m-1}) \\
	&\quad\quad- ( t_{n-m}^{2} + t_{n-m-1}^{2}) + 2 (\Delta t)^{2} ) 
	\chi_{E_{n-m-1}}(x,y)\quad=: \Upsilon^{n-m}(x,y)\ .
\end{align*}
Therefore we get:
\begin{equation}\label{eq:B}
	B = \sum_{m=1}^{N_t} \sum_{i=1}^{N_s} \frac{c_{i}^{m} }{2 \pi} \int_{\Gamma \times \Gamma} 
	\frac{1}{|x-y|} \operatorname{curl}_{|\Gamma} \varphi^1_{i}(y) \operatorname{curl}_{|\Gamma} \varphi^1_{j}(x)\ \Upsilon^{n-m}(x,y)\ ds_y\ ds_x.
\end{equation}
In summary, both terms \( A \) and \( B \) can be expressed as integrals over the light cones \( E_{n-m} \), \( E_{n-m-1} \), and \( E_{n-m-2} \), associated with the \( m \)-th ansatz and the \( n \)-th test function in time.
Due to causality, the terms in \( A \) and \( B \) vanish for \( t_n > t_m \). \\
{As described for the adjoint double-layer operator $\mathcal{K'}$ above, the accurate numerical evaluation of $A$ and $B$ is based on a geometric decomposition of the light cone combined with a composite $hp$ quadrature rule \cite{ostermann2010numerical,iccs2008,gimperlein2016adaptive}.}

\bibliographystyle{plain}

\bibliography{lit3}

\end{document}